\documentclass[11pt]{amsart}
\usepackage[hmargin=32mm, vmargin=27mm, includefoot, 
]{geometry}
\usepackage[
colorlinks,bookmarksopen=true]{hyperref}
\usepackage[latin1]{inputenc}  
\usepackage[T1]{fontenc}       
\usepackage{amssymb}
\usepackage{latexsym}
\usepackage{mathrsfs}
\usepackage{xspace}
\usepackage{graphics, graphicx, epsfig}
\usepackage[all, cmtip]{xy}
\newtheorem{prop}{Proposition}[section]
\newtheorem{thm}[prop]{Theorem}
\newtheorem*{thm*}{Theorem}

\newtheorem*{addendum*}{Addendum}
\newtheorem{cor}[prop]{Corollary}
\newtheorem{lem}[prop]{Lemma}

\newtheorem*{convention*}{Convention}
\theoremstyle{definition}
\newtheorem*{defn*}{Definition}

\newtheorem{remark}[prop]{Remark}

\newtheorem*{scholium*}{Scholium}
\theoremstyle{remark}

\newtheorem*{example*}{Example}
\numberwithin{equation}{section}

\newcommand{\FF}{\mathbf{F}}
\newcommand{\GG}{\mathbf{G}}

\newcommand{\HH}{\mathbf{H}}

\newcommand{\NN}{\mathbf{N}}
\newcommand{\OO}{\mathbf{O}}

\newcommand{\QQ}{\mathbf{Q}}
\newcommand{\RR}{\mathbf{R}}

\newcommand{\ZZ}{\mathbf{Z}}

\newcommand{\inv}{^{-1}}

\newcommand{\norma}{\mathscr{N}}

\newcommand{\Comm}{\mathrm{Comm}}

\newcommand{\se}{\subseteq}

\def\bs#1.{
              \def\temp{#1}
              \ifx\temp\empty
                   \mathcal{B}
              \else
                   \mathcal{B}(#1)
              \fi
}
\newcommand{\cat}{{\upshape CAT(0)}\xspace}
\newcommand{\tangle}[2]
{\angle_\mathrm{T}(#1,#2)}
\newcommand{\aangle}[3]
{\angle_{#1}(#2,#3)}
\newcommand{\cangle}[3]
{\overline{\angle}_{#1}(#2,#3)}
\DeclareMathOperator{\chr}{char}
\DeclareMathOperator{\rank}{rank}
 
\DeclareMathOperator{\Stab}{Stab}

\DeclareMathOperator{\Isom}{Is}

\DeclareMathOperator{\charact}{char}
\newcommand{\bd}{\partial}

\def\Aut{\mathop{\mathrm{Aut}}\nolimits}

\begin{document}
\title[Lattices in products of Kac--Moody groups]{A lattice in more than two Kac--Moody groups is arithmetic}
\author[Pierre-Emmanuel Caprace]{Pierre-Emmanuel Caprace*}
\address{UCLouvain, 1348 Louvain-la-Neuve, Belgium}
\email{pierre-emmanuel.caprace@uclouvain.be}
\thanks{*F.N.R.S. Research Associate}
\author[Nicolas Monod]{Nicolas Monod$^\ddagger$}
\address{EPFL, 1015 Lausanne, Switzerland}
\email{nicolas.monod@epfl.ch}
\thanks{$^\ddagger$Supported in part by the Swiss National Science Foundation}
\begin{abstract}
Let $\Gamma < G_1 \times \dots \times G_n$ be an
irreducible lattice in a product of infinite
irreducible complete Kac--Moody groups of simply laced type over finite fields. We show that if $n \geq 3$, then
each $G_i$ is a simple algebraic group over a local field and $\Gamma$ is an $S$-arithmetic lattice. This relies on
the following alternative which is satisfied by any
irreducible lattice provided $n \geq 2$:
either $\Gamma$ is an $S$-arithmetic (hence linear) group, or $\Gamma$ is not residually finite. In that case, it is
even virtually simple when the ground field is large enough.

More general \cat groups are also considered throughout.
\end{abstract}
\maketitle
\let\languagename\relax  


\section{Introduction}\label{sec:intro}

The theory of lattices in semi-simple Lie and algebraic groups has witnessed tremendous developments over the past
fourty years. It has now reached a remarkably deep and rich status, notably thanks to the celebrated work of
G.~Margulis, whose main aspects may be consulted in~\cite{Margulis}. Amongst the followers and exegetes of
Margulis' work, several authors extended the methods and results pertaining to this classical setting to broader
classes of lattices in locally compact groups.
It should be noted however that as of today there
exists apparently no characterisation of the $S$-arithmetic lattices purely within the category of lattices in compactly
generated locally compact groups.

\medskip

It turns out that relatively few examples of compactly generated topologically simple groups are known to
possess lattices; to the best of our knowledge, they are all \emph{locally compact \cat groups}.
In fact, the only examples which are neither algebraic nor Gromov hyperbolic are all
automorphism groups of non-Euclidean locally finite buildings. Amongst these, the most prominent family consists
perhaps of the so-called \textbf{irreducible complete Kac--Moody groups} over finite fields constructed by
J.~Tits~\cite{Tits87} (see \S\,\ref{sec:KM} below for more details).

\bigskip

We now proceed to describe our main result. To this end, fix a positive integer $n$.

For each $i \in \{1, \dots, n\}$, let $X_i$ be a proper \cat space and $G_i < \Isom(X_i)$ be a closed subgroup acting cocompactly.

\begin{thm}\label{thm:RankRigid}
Let  $\Gamma < G_1 \times \dots \times G_n$ be any lattice whose projection to each $G_i$ is faithful. Assume
that $G_1$ is an irreducible complete Kac--Moody group of simply laced type over a finite field.

If $n \geq 3$, then each $G_i$ contains a cocompact normal subgroup which is a compact
extension of a semi-simple group over a local field, and $\Gamma$ is an $S$-arithmetic lattice.

\smallskip
(The same conclusion holds if $G_1$ is instead assumed to be a non-Gromov-hyperbolic irreducible complete Kac--Moody group
of $3$-spherical type over a finite field of characteristic~$\neq 2$.)
\end{thm}

One can summarise the above result as follows: \itshape As soon as $n\geq 3$ and one of the factors $G_i$ is Kac--Moody as above,
all $G_i$ are topologically commensurable to semi-simple algebraic groups and the lattice is $S$-arithmetic.\upshape

\begin{remark}
Our aim in Theorem~\ref{thm:RankRigid} is to provide a statement without any restrictive assumptions on the lattice
$\Gamma$, on the spaces $X_i$ or on the groups $G_i$. Considerable complications are caused by the fact that $\Gamma$
is not supposed finitely generated. It turns out \emph{a posteriori} that only finitely generated lattices exist~---
as a consequence of arithmeticity.
\end{remark}

\begin{remark}
The assumption on faithfulness of the projections of $\Gamma$ to each factor $G_i$ is a form of irreducibility.
We refer to Section~\ref{sec:irred} below for a  discussion of the different possible definitions of irreducibility for lattices in
products of locally compact groups.
\end{remark}

The above theorem is a new manifestation of the phenomenon that ``high rank''
yields rigidity. Numerous other results support this vague statement, including the rank rigidity of Hadamard
manifolds, the arithmeticity of lattices in higher rank semi-simple groups, or the fact that any irreducible
spherical building of rank~$\geq 3$ (resp. affine building of dimension~$\geq 3$) is associated to a simple
algebraic group (possibly over a skew field).

\medskip

Theorem~\ref{thm:RankRigid} will be established with the help of the following \textbf{arithmeticity vs. non-residual-finiteness
alternative}.

\begin{thm}\label{thm:alternative:1}
Let $\Gamma < G_1 \times \dots \times G_n$ be a lattice which is algebraically irreducible. Assume that
$G_1$ is an infinite irreducible complete Kac--Moody group of simply-laced type over a finite field.

If $n \geq 2$ then either $\Gamma$ is an $S$-arithmetic group or $\Gamma$ is not residually finite.

\smallskip
(The same conclusion holds if $G_1$ is instead assumed to be a non-Gromov-hyperbolic irreducible complete Kac--Moody group
of $3$-spherical type over a finite field of characteristic~$\neq 2$.)
\end{thm}

It is known that if $G = G_1 \times G_2$ is a product of two isomorphic complete Kac--Moody groups over a
sufficiently large finite field, then $G$ contains at least one irreducible non-uniform lattice (see
\cite{RemCRAS},~\cite{CarGar}). In~\cite{CaRe}, this specific example is shown to be simple provided $G_1$ and
$G_2$ are non-affine (and without any other restriction on the type). Theorem~\ref{thm:alternative:1} shows in
particular that, under appropriate assumptions, \emph{all} irreducible lattices in $G$ are virtually simple.
More precisely, we have the following \textbf{arithmeticity vs. simplicity alternative} which, under more precise
hypotheses, strengthens the alternative from Theorem~\ref{thm:alternative:1}.

\begin{cor}\label{cor:simple}
Let $G= G_1 \times \dots \times G_n$, where $G_i$ is an infinite irreducible Kac--Moody group of
simply-laced type over a finite field $\FF_{q_i}$  (or a non-affine non-Gromov-hyperbolic irreducible complete Kac--Moody
group of $3$-spherical type over a finite field $\FF_{q_i}$ of characteristic~$\neq 2$). Let $\Gamma < G$
be a topologically irreducible lattice; if $\Gamma$ is not uniform, assume in addition that $q_i \geq 1764^{d_i}/25$
for each $i$, where $d_i$ denotes the maximal rank of a finite Coxeter subgroup of
the Weyl group of $G_i$. If $n \geq 2$, then exactly one of the following assertions holds:
\begin{enumerate}
\item Each $G_i$ is of affine type and $\Gamma$ is an arithmetic lattice.

\item $n=2$ and $\Gamma$ is virtually simple.
\end{enumerate}
\end{cor}

It is important to remark that, in most cases, a group $G$ satisfying the hypotheses of any of the above statements
does \emph{not} admit any uniform lattice (see Remark~\ref{rem:BorelHarder} below).

\begin{cor}\label{cor:MaximalLattices}
Let $G= G_1 \times \dots \times G_n$, where $G_i$ is an infinite irreducible non-affine Kac--Moody groups of
simply-laced type over a finite field $\FF_{q_i}$ (or a non-affine non-Gromov-hyperbolic irreducible complete Kac--Moody
group of $3$-spherical type over a finite field $\FF_{q_i}$ of characteristic~$\neq 2$). Assume that $q_i \geq
1764^{d_i}/25$ for each $i$, where $d_i$ denotes the maximal rank of a finite Coxeter subgroup of the
Weyl group of $G_i$.

If $n\geq 2$, then any topologically irreducible lattice of $G$ has a discrete commensurator, and is thus contained in a unique maximal
lattice.
\end{cor}

\medskip
Our proof of Theorems~\ref{thm:RankRigid} and~\ref{thm:alternative:1} builds upon the general methods developed
in~\cite{CaMoa,CaMob} for studying lattices in isometry groups of non-positively curved spaces. Our treatment of the  residual-finiteness/simplicity
dichotomy is inspired by the work of Burger--Mozes for tree lattices~\cite{Burger-Mozes2}.

\subsection*{Acknowledgement} We are grateful to the anonymous referee for his/her useful comments.

\section{Lattices in non-positive curvature}

\subsection{The set-up}\label{sec:set-up}%
We now introduce the setting for this section and the subsequent ones. The situation will differ from the very
simple assumptions made in the Introduction; indeed our first task in the proof of Theorems~\ref{thm:RankRigid}
and~\ref{thm:alternative:1} will be to reduce them to the set-up below.

Fix an integer $n\geq 2$. For
each $i \in \{1, \dots, n\}$, let $X_i$ be an irreducible proper \cat space not isometric to the real line;
irreducibility of $X_i$ means that it does not split (non-trivially) as a direct product. It follows that $X_i$
has no Euclidean factor. We also assume that the boundary $\bd X_i$ is finite-dimensional for the Tits topology
(though this assumption will only be used in later sections).

Let further $G_i < \Isom(X_i)$ be a compactly generated closed subgroup without global fixed point at infinity. We assume that $G_i$
acts \textbf{minimally} in the sense that there is no invariant closed convex proper subspace of $X_i$. We point
out that this assumption is automatically fulfilled upon passing to some subspace since there is no fixed point
at infinity, compare Remark~36 in~\cite{Monod_superrigid}.

We set $G = G_1 \times \dots \times G_n$ and $X = X_1 \times \dots \times X_n$. Finally,
let $\Gamma < G$ be a lattice.

\medskip

It was established in~\cite{CaMob} that the following ``Borel density'' holds.

\begin{prop}\label{prop:LatticeConsequences}
The action of $\Gamma$ and its finite index subgroups on $X$ is minimal and without fixed points at infinity.
\end{prop}

\begin{proof}
This is a special case of Theorem~2.4 in~\cite{CaMob}.
\end{proof}

\subsection{Irreducible lattices and residual finiteness}\label{sec:irred}

Let $G = G_1 \times \cdots \times G_n$ be a locally compact group. The following properties provide several
possible definitions of irreducibility for a lattice $\Gamma$ in the product $G= G_1 \times \cdots \times G_n$, which
we shall subsequently discuss.

\begin{itemize}
\item[\textbf{(Irr1)}] The projection of $\Gamma$ to any proper sub-product of $G$ is dense (and all $G_i$ are non-discrete).
In this case $\Gamma$ is called \textbf{topologically irreducible}.

\item[\textbf{(Irr2)}] The projection of $\Gamma$ to each factor $G_i$ is injective.

\item[\textbf{(Irr3)}] The projection of $\Gamma$ to any proper sub-product of $G$ is non-discrete.

\item[\textbf{(Irr4)}] $\Gamma$ has no finite index subgroup which splits as a direct product of two infinite
subgroups. In this case $\Gamma$ is called \textbf{algebraically irreducible}.
\end{itemize}

It turns out, as is well known, that if each factor $G_i$ is a semi-simple linear algebraic group, then all four
properties (Irr1)--(Irr4) are equivalent. We shall show that, in the setting of
\S\,\ref{sec:set-up}, the implications (Irr2)$\Rightarrow$(Irr3)$\Rightarrow$(Irr4) do hold. If one assumes
furthermore that each $G_i$ is topologically simple and compactly generated, then each $G_i$ has trivial
quasi-centre by~\cite[Theorem~4.8]{BarneaErshovWeigel} and, hence, one has (Irr1)$\Rightarrow$(Irr2) in that
case. However,   even in the setting of \S\,\ref{sec:set-up}, \emph{none} of the implications
(Irr2)$\Rightarrow$(Irr1), (Irr3)$\Rightarrow$(Irr2) or (Irr4)$\Rightarrow$(Irr3) holds true.

A crucial point in the proof of Theorem~\ref{thm:alternative:1} is that, nevertheless, the implication
(Irr4)$\Rightarrow$(Irr2) becomes true provided the lattice $\Gamma$ is residually finite, see
Proposition~\ref{prop:ResFinite:irred} below.

\medskip

From now on, we retain the notation of \S\,\ref{sec:set-up}.

The following result implies that (Irr2)$\Rightarrow$(Irr3)$\Rightarrow$(Irr4). 

\begin{prop}\label{prop:irred}
\ 
\begin{enumerate}
\item[(i)] If the projection of $\Gamma$ to each factor $G_i $ is faithful, then the projection of $\Gamma$ to any proper sub-product of $G$ is non-discrete.

\item[(ii)] If the projection of $\Gamma$ to any proper sub-product of $G$ is non-discrete, or if the projection to at least one factor $G_i$ is faithful, then $\Gamma$ is algebraically irreducible.
\end{enumerate}
\end{prop}

\begin{proof}
(i) Assume that the projection of $\Gamma$ to each factor $G_i $ is faithful and consider any (non-trivial) regrouping of factors $G=G'\times G''$; we need to show that the projection of $\Gamma$ to $G'$ is not discrete. Assume thus that the latter is discrete. In that case, Lemma~I.1.7 in~\cite{Raghunathan} ensures that $\Gamma \cap ( \{1\} \times G'')$ is a lattice in $\{1\} \times G''$. In particular it is non-trivial. Therefore the projection of $\Gamma$ to $G'$ cannot be faithful (a fortiori to some $G_i$).

\bigskip \noindent 
(ii) Suppose that some finite index subgroup $\Gamma_0<\Gamma$ admits a splitting.
Proposition~\ref{prop:LatticeConsequences} implies in particular that $\Gamma_0$ acts on $X$ as well as on each factor $X_j$, minimally
and without fixed points at infinity. These are exactly the assumptions necessary to apply the splitting
theorem of~\cite{Monod_superrigid}. The latter provides a splitting of $X$ as $X = Y \times Z$ which, by the canonicity of the geometric decomposition  $X \cong X_1 \times \dots \times X_n$, must correspond to some regrouping of irreducible factors of $X$. In other words we have a non-empty subset $J \subsetneq \{1, \dots, n\}$ such that $Y = \prod_{j \in J} X_j $ and $Z = \prod_{j \not \in J} X_j$.  The desired result now follows from the fact that  the respective $\Gamma_0$-actions on $Y$ and $Z$ are discrete but not faithful. 
\end{proof}

Examples showing that the implication (Irr3)$\Rightarrow$(Irr2) fails in the setting of \S\,\ref{sec:set-up} can be obtained as extensions of arithmetic lattices by products or free groups, using similar constructions as in \cite[\S6.C]{CaMob} (suggested by Burger--Mozes). The following result  shows that (Irr4)$\Rightarrow$(Irr3) provided that the lattice $\Gamma$ is finitely generated. 

\begin{prop}\label{prop:NonDiscreteProj}
Assume that $\Gamma$ is finitely generated and algebraically irreducible. Then the projection of $\Gamma$ to any proper sub-product of
$G = G_1 \times \cdots \times G_n$ has non-discrete image. 
\end{prop}
\begin{proof}
See Theorem~4.2(i) in~\cite{CaMob}.
\end{proof}

We shall now describe an example showing that the implication (Irr4)$\Rightarrow$(Irr3) fails to hold if one removes the hypothesis that $\Gamma$ be finitely generated. This also illustrates some of the technical difficulties that are unavoidable in the proof of our main results, since we deal with general (\emph{i.e.} possibly non-uniform infinitely generated) lattices.

\begin{example*}
Let $A = \bigoplus_{n\geq 0} \ZZ/2$ and $M = A * A$. Then $M$ can be realised as a non-uniform lattice in the group $\Aut(T_3)$ of the regular tree of degree~$3$. To see this, one can express $M$ as the fundamental group of a graph of groups as follows. Let $\lambda^+$ and $\lambda^-$ be two copies of the simplicial half-line, and let us index the respective vertices of $\lambda^+$ and $\lambda^-$ by the strictly positive integers in the natural order. On both $\lambda^+$ and $\lambda^-$, we attach the group $(\ZZ/2)^n$ to the vertex $n$ and to the edge joining $n$ to $n+1$. The embedding of the edge group attached with $[n, n+1]$ to the vertex group attached with $n+1$ is the natural inclusion of $(\ZZ/2)^n$ in $(\ZZ/2)^{n+1}$ defined by $x \mapsto (x, 0)$. Finally, we join the vertex $1$ of $\lambda^+$ to the vertex $1$ of $\lambda^-$ by an edge to which we attach the trivial group. In this way, we obtain a graph of groups, which is simplicially isomorphic to the line. Its 
 fundamental group is isomorphic to $M$ and acts on the universal cover $T_3$ as a non-uniform lattice. 

Similarly, we set $B = \bigoplus_{n \geq 0} \ZZ/3$ and view the group $ N = B*B$ acting as a non-uniform lattice on the regular tree $T_4$ of degree $4$ using the same construction, but replacing $\ZZ/2$ by $\ZZ/3$. 

Now we define an action of $M$ by automorphisms on $N$. Clearly $A$ acts on $B$ by non-trivial automorphisms componentwise, so that the semi-direct product $A \ltimes B$ is isomorphic to $\bigoplus_{n \geq 0} (\ZZ/2 \ltimes \ZZ/3)$: in each coordinate, the group $\ZZ/2$ acts on $\ZZ/3$ as the (unique) non-trivial automorphism. This action extends naturally to a diagonal action of $A$ on $B \times B$ which, post-composed with the embedding of sets $B \times B \hookrightarrow B*B$, defines an action of $A$ on the generators of $N  = B*B$ preserving all the defining relations. Thus $A$ acts on $N = B*B$ by automorphisms. Precomposing this with the natural quotient map $M = A*A \to A$ which annihilates the second free factor, we obtain a homomorphism
$$
\alpha : M \to \Aut(N).
$$
Since the $M$-action on $N$ preserves the graph of group decomposition of $N$, it extends to an $M$-action by automorphisms on $T_4$ which, by a slight abuse of notation, we also denote by $\alpha$. As a subgroup of $\Aut(T_4)$, the group $\alpha(M)$ fixes pointwise a line; the closure of $M$ in $\Aut(T_4)$ is a compact subgroup $Q$ isomorphic to $\prod_{\ZZ} \ZZ/2$. 

Set now 
$$
 \Gamma = N \rtimes_\alpha M
  \hspace{1cm} \text{and} \hspace{1cm}
  G = \Aut(T_4) \times \Aut(T_3).$$
We have already defined an embedding $f_4 : \Gamma \to \Aut(T_4)$ and a lattice embedding of $M$ in $\Aut(T_3)$. Precomposing the latter with the quotient map $\Gamma \to M$,  we obtain a homomorphism $f_3 :  \Gamma \to \Aut(T_3)$ whose image is the lattice $M < \Aut(T_3)$. Finally, we define an injective homomorphism
$$
f : \Gamma \to G : \gamma \mapsto (f_4(\gamma), f_3(\gamma)).
$$
The image $f(\Gamma)$ is discrete. Moreover, since the image of $f(\Gamma)$ is a lattice in $\Aut(T_3)$ and the kernel of the projection of $f(\Gamma)$ to $\Aut(T_3)$ is a lattice in $\Aut(T_4)$, it follows that $f(\Gamma)$ is a lattice in $G$. 

Remark that $\Gamma$ is algebraically irreducible since no finite index subgroup of $M$ is normal in $\Gamma$. The projection of $\Gamma$ to $\Aut(T_3)$ is discrete while its projection to $\Aut(T_4)$ is not, since its closure is isomorphic to $N \rtimes Q$. This shows that Proposition~\ref{prop:NonDiscreteProj} does not hold for lattices which are not finitely generated.
\end{example*}

\bigskip
We finish this subsection with a crucial ingredient in the proof of Theorem~\ref{thm:alternative:1} which shows that the implication
(Irr4)$\Rightarrow$(Irr2) does however hold under the extra assumption that the lattice $\Gamma$ is residually finite~--- even if it is not finitely generated. 

\begin{prop}\label{prop:ResFinite:irred}
Assume that $\Gamma$ is residually finite and algebraically irreducible. Then the projection of $\Gamma$ to each
$G_i $ is faithful.
\end{prop}

\begin{proof}
In the special case when $\Gamma$ is finitely generated, we obtained this result in Theorem~4.10 of~\cite{CaMob}.
In the present level of generality, we can write $\Gamma$ as the union of an ascending sequence of finitely generated subgroups
$(\Gamma_j)_{j \geq 0}$ because $\Gamma$ is countable since $\Isom(X)$ is second countable.

We let $H_i$ denote the closure of the projection of $\Gamma$ to $G_i$. 
Upon reordering the factors, we may assume that there is some $s \in \{0, \dots, n\}$ such that $H_i$ is discrete if and only if $i > s$. We remark that if $s=0$, then $H_i $ is discrete for all $i$. Therefore, Lemma~I.1.6 from \cite{Raghunathan} ensures that  $H_1 \times \dots \times H_n$ is a lattice in $G=G_1 \times \dots \times G_n$ and that the index of $\Gamma$ in $H_1 \times \dots \times H_n$ is finite. Therefore the product group $(\Gamma \cap H_1) \times \dots \times (\Gamma \cap H_n)$ has finite index in $\Gamma$, contradicting the fact that $\Gamma $ is algebraically irreducible. Thus $s >0$ as asserted. 

By~\cite[Corollary~1.11]{CaMoa}, each $G_i$ is either totally disconnected or an adjoint simple non-compact Lie group. By the definition of $s$, the group $H_i$ is non-discrete for all $i \leq s$ and, hence, dense in every connected factor of $ H_1 \times \dots \times H_s$ by Borel density~\cite[4.2]{Borel60} (see also~\cite[II.6.2]{Margulis}). Thus, for all $i \leq s$, the group $H_i$ is a non-discrete closed subgroup which coincides with $G_i$ if the latter is not totally disconnected.

\medskip
Let $I \subseteq \{1, \dots, s\}$ be any non-empty subset. We claim that if the projection of $\Gamma$ to $\prod_{i \not \in I} H_i$ is not faithful, then  the projection of $\Gamma$ to $\prod_{i \in I} H_i$ is discrete. 

\medskip
In order to establish the claim, we let $C$ denote the closure of the projection of $\Gamma$ to $\prod_{i  \in I} H_i$ and let 
$$
N = \Gamma \cap (\prod_{i  \in I} H_i \times \prod_{i \not \in I} \{1\}) < H_1 \times \dots \times H_n.
$$
Then $C$ is totally disconnected; this is shown in the proof of
Theorem~4.10 in~\cite{CaMob} by arguments that do not depend on the finite generation of $\Gamma$, but use the fact that $H_i$ is either totally disconnected or a connected simple Lie group for all $i \in I$.

We now assume that $N$ is non-trivial and need to deduce that $C$ is discrete. Since $\Gamma$ has trivial amenable radical~\cite[Corollary~2.7]{CaMob} the normal subgroup $N \unlhd
\Gamma$ is not locally finite and, hence, we can assume upon discarding the first few indices in the filtration $(\Gamma_j)_{j \geq 0}$ that
$\Gamma_j \cap N$ is infinite for each $j\geq 0$. Furthermore, since $\Gamma$ has no fixed point $\bd X$ by
Proposition~\ref{prop:LatticeConsequences} and since $\bd X$ is compact when endowed with the cone topology, we can moreover assume that $\Gamma_j$ has no fixed point in $\bd X$ for all $j$. Finally, let $Q <
C$ be a compact open subgroup and denote by $C_j<C$ the subgroup generated by $Q$ and the image of $\Gamma_j$.

By construction the group $C_j$ is compactly generated, it acts without fixed point at infinity and the
intersection of $C_j$ with the image of $N$ in $C$ is an infinite discrete normal subgroup of $C_j$. Since $\Gamma \cap \big((\prod_{i \not \in I} H_i) \cdot C_j\big)$ projects
densely to $C_j$, we deduce from~\cite[Proposition~4.8]{CaMob}
that  $[C_j \cap N, C_j^{(\infty)}]=1$, where $C_j^{(\infty)}$ denotes the intersection of
all open normal subgroups of $C_j$.
We recall that a non-compact group of isometries of a proper \cat space $X$
acting without fixed point at infinity has a compact centraliser in $\Isom(X)$ (though this is an overkill, it follows
\emph{e.g.} from the splitting theorem in~\cite{Monod_superrigid}); hence $C_j^{(\infty)}$ is compact.

On the other hand, the group $C_j$ possesses a maximal compact normal subgroup because it acts without fixed
point at infinity; this follows \emph{e.g.} from Corollary~5.8 in~\cite{CaMoa}, the compact subgroup being
the kernel of the $C_j$-action on a minimal subspace.
Therefore, it follows from~\cite[Proposition~6.12]{CaMoa} that $C_j^{(\infty)}$ is non-compact
whenever $C_j$ is non-discrete. This shows that $C_j$ is discrete. By construction $C_j$ is open in $C$, thus
$C$ is discrete as well. This proves the claim. 

\medskip
We shall now establish that $s=n$. To this end, we notice that the projection of $\Gamma$ to $G_{s+1}\times\cdots\times G_n$ cannot be faithful since it has discrete image (see \cite[Lemma~I.1.7]{Raghunathan}). Applying the claim with the set $I = \{1, \dots, s\}$, we infer that the closure $C$ of the projection of $\Gamma$ to $ H_1 \times \dots \times H_s$ is discrete. By the definition of $s$, we also observe that the closure $D$ of the projection of $\Gamma$ to $H_{s+1} \times \dots \times H_n$ is discrete. Thus $\Gamma$ is contained in the discrete subgroup $C \times D < (H_1 \times \dots \times H_s) \times (H_{s+1} \times \dots \times H_n)$. Since $\Gamma $ is a lattice, it must therefore have finite index in $C \times D$. Therefore $\Gamma$ has a finite index subgroup which splits as a direct product, which contradicts the hypothesis that $\Gamma$ is algebraically irreducible. This contradiction confirms that $s=n$. 

\medskip
Finally, assume for a contradiction that the projection of $\Gamma $ to $G_k$ is not faithful for some $k \in \{1, \dots, n\}$. We then invoke the claim above to the set $I =\{1, \dots, n\} \setminus \{k\}$. From the claim, we infer that the projection $C'$ of $\Gamma$ to $\prod_{i \neq k} H_i$ is discrete. Therefore, this projection cannot be faithful (see \cite[Lemma~I.1.7]{Raghunathan}). We can then invoke the claim one more time, now with the set $I = \{k\}$. This implies that the projection $D'$ of $\Gamma$ to $H_k$ is discrete, contradicting $s=n$.
\end{proof}

\subsection{Open subgroups}\label{sec:open}
Recall that, given a lattice $\Lambda$ in a locally compact group $H$ and any open subgroup $P < H$, the
intersection $\Lambda \cap P$ is a lattice in $P$; indeed a Haar measure for $P$ may be obtained by restricting the Haar measure of $H$.
Furthermore, if $\Lambda$ is uniform in $H$, so is $\Lambda \cap P$ in $P$. 
We shall frequently take advantage of this basic observation and study the intersection $\Gamma \cap P$ for various open
subgroups $P < G$.

\begin{lem}\label{lem:intermediate:lattice}
Let $H, J$ be locally compact groups, $\Lambda<H\times J$ a lattice, $P<H$ an open subgroup and $\Lambda_P=\Lambda\cap (P\times J)$.
Then any intermediate group $\Lambda_P < \Lambda' < \Lambda$ is a lattice in $H'\times J$ and in $H'\times J'$, where $H'$ and $J'$ are
the closure of the projection of $\Lambda'$ to $H$ and $J$ respectively.
\end{lem}

\begin{proof}
Let $H^*$ be the the closure of the projection of $\Lambda$ to $H$ and $P^*=P\cap H^*$. Then $\Lambda$ is a lattice in $H^*\times J$
by~\cite[I.1.6]{Raghunathan}. Moreover, $\Lambda_P$ is a lattice in $P^*\times J$ projecting densely to $P^*$ since $P$ is open;
in particular, $P^*\se H'$. Let $F\se P^*\times J$ be a (left) fundamental domain for $\Lambda_P$ in $P^*\times J$.

We claim that the $\Lambda'$-translates of $F$ cover
$H'\times J$. Pick thus any $(h_0, j_0)$ in $H'\times J$. Since $P^*$ is open in $H'$, there is $(h_1, j_1)$ in $\Lambda'$ such that
$h_1 h_0\in P^*$. Since $(h_1 h_0, j_1 j_0)\in P^*\times J$, there is $(h_2, j_2)$ in $\Lambda_P$ such that
$(h_2 h_1 h_0, j_2 j_1 j_0)\in F$. Since $(h_2 h_1, j_2 j_1)\in \Lambda'$, this proves the claim.

Since $\Lambda'$ is discrete and since the Haar measures of $P^*\times J$ extend to Haar measures of $H'\times J$,
we conclude that $\Lambda'$ is indeed a lattice in $H'\times J$. Applying again~\cite[I.1.6]{Raghunathan},
we deduce that it is also a lattice in $H'\times J'$.
\end{proof}

We now return to our geometric setting.

\begin{prop}\label{prop:CompactOpen:irred}
Let $P < G_1$ be an open subgroup and set
$$ \Gamma_{P} = \Gamma \cap (P \times G_2 \times \dots \times G_n).$$
Assume that the projection of $\Gamma_P$ to some $G_i$ with $i\geq 2$ is faithful.
Then $\Gamma_{P}$ is algebraically irreducible.
\end{prop}

\begin{proof}
In order to argue as in Proposition~\ref{prop:irred}, we need to show that the $\Gamma_P$-action
on $X_i$ is minimal and without fixed point at infinity.

We claim that without loss of generality we may assume that $G_1$ is totally disconnected. Indeed, otherwise by
Corollary~1.11 in~\cite{CaMoa} the group $G_1$ is an almost connected simple Lie group. In that case the open
subgroup $P < G_1$ has finite index in $G_1$ and hence $\Gamma_{P}$ has finite index in $\Gamma$. The claimed
statement is thus a case of Proposition~\ref{prop:irred}.

In view of the claim, we assume that $G_1$ is totally disconnected; hence so is $P$.
Therefore there exists a compact open subgroup $U<P$ (see~\cite[III.4 No~6]{BourbakiTGI}).
Then the group
$$ \Gamma_{U} = \Gamma \cap (U \times G_2 \times \dots \times G_n)$$
is a lattice in $U \times G_2 \times \dots \times G_n$ and thus its projection to $H := G_2 \times \dots \times G_n$
is a lattice as well. Since the $H$-action on $Y := X_2 \times \cdots \times X_n$ is minimal and without fixed point
at infinity, so is the $\Gamma_{U}$-action by Proposition~\ref{prop:LatticeConsequences}. Now we deduce \emph{a fortiori} that
the $\Gamma_{P}$-action on $Y$ and hence on $X_i$ is minimal and without fixed point at infinity.
\end{proof}

\begin{cor}\label{cor:Open:irred}
Let $P < G_1$ be any open subgroup and set
$$ \Gamma_{P} = \Gamma \cap (P \times G_2 \times \dots \times G_n).$$
If $\Gamma$ is algebraically irreducible and residually finite, then so is $\Gamma_{P}$.

\end{cor}

\begin{proof}
By Proposition~\ref{prop:ResFinite:irred} the projection of $\Gamma$ to each $G_i$ is faithful. Thus we may
apply Proposition~\ref{prop:CompactOpen:irred}.
\end{proof}

\subsection{Cofinite embeddings of semi-simple groups}
We do not know if a semi-simple algebraic group can appear as a subgroup of finite covolume in a locally compact
group without being cocompact\footnote{Note added in proof: we have been informed that Bader--Furman--Sauer address this question in forthcoming work.}. We shall prove that this does not happen in the \cat setting.

\begin{prop}\label{prop:semi-simple:cocompact}
Let $H$ be a locally compact group and $L<H$ a closed subgroup of finite covolume which is a compact extension
of a semi-simple algebraic group. Suppose that $H$ admits a cocompact proper continuous isometric action on some
\cat space.

Then $H/L$ is compact.

Morover, if the semi-simple group has no rank one factors, then upon factoring out a (unique maximal) compact
normal subgroup, $H$ is a group of automorphismsm of the semi-simple algebraic group.
\end{prop}

The following fact is well-known.

\begin{lem}\label{lem:proba}
A group of isometries preserving a non-zero finite measure on a complete \cat space fixes a point.
\end{lem}

\begin{proof}
Let $G$ be the group, $X$ the space and $\mu$ the measure. There is a bounded set $B\se X$ such that
$\mu(B)>\mu(X)/2$. Therefore $gB\cap B \neq \varnothing$ for all $g\in G$. It follows that $G$ has a bounded
orbit and hence a fixed point by Cartan's circumcentre principle~\cite[II.2.8]{Bridson-Haefliger}.
\end{proof}

\begin{lem}\label{lem:NoEuclidean}
Let $H$ be a locally compact group containing a closed subgroup of finite covolume which is a compact extension
of a semi-simple algebraic group. Then any continuous isometric $H$-action on a proper \cat space preserves a
non-empty closed convex subset with trivial Euclidean factor.
\end{lem}

\begin{proof}
Invoking repeteadly the canonical Euclidean decomposition~\cite[II.6.15]{Bridson-Haefliger}, it suffices to
prove that any continuous isometric $H$-action on any $\RR^d$ has a fixed point. Let $L<H$ be the given subgroup
with $K_0\lhd L$ compact normal such that $L/K_0$ is semi-simple. The non-empty subspace of $K_0$-fixed points
is an affine subspace preserved by $L$; we therefore have an isometric action of the semi-simple group $L/K_0$
on some $\RR^{d'}$.

It is well-known that all such $L/K_0$-actions have a fixed point. Therefore $L$ fixes a point in
$\RR^d$, which implies by Lemma~\ref{lem:proba} that $H$ also fixes a point.
\end{proof}

\begin{proof}[Proof of Proposition~\ref{prop:semi-simple:cocompact}]
Let $X$ be a \cat space as in the statement; it is necessarily proper. Since $H$ acts cocompactly, it has a
minimal convex invariant subspace and thus we can assume $X$ minimal upon factoring out the compact kernel of
the $H$-action on that subspace. We note in passing that this kernel is a (necessarily unique) maximal compact
normal subgroup of $H$.

Lemma~\ref{lem:NoEuclidean} implies that $X$ has trivial Euclidean factor. Moreover, we claim that $H$ has no fixed point at infinity.
Indeed, otherwise by minimality the corresponding Busemann character $H\to \RR$ would be non-trivial. This however would produce
a non-trivial character of $L$  which would thus descend non-trivially to the semi-simple group, which is absurd.

By Theorem~2.4 in~\cite{CaMob}, the $L$-action on $X$ is minimal and without fixed point at infinity. In particular, $L$ has no non-trivial compact
normal subgroup and we can decompose it into its simple factors $L=L_1\times \cdots\times L_n$. Each $L_i$ is
non-compact and we can assume $n\geq 1$ since otherwise $H$ is compact (actually trivial at this point).

In view of Addendum~1.8 in~\cite{CaMoa} we can write $X=X_1\times \cdots\times X_n$, where each $L_i$ acts
minimally on $X_i$; the finite-dimensionality of $\bd X$ holds by Theorem~C in~\cite{Kleiner} since $H$ acts cocompactly.
Moreover, each $\bd X_i$ is finite-dimensional and each $L_i$ has full limit set in $\bd X_i$ because the two corresponding statements for
$\bd X$ and the $L$-action on $X$ hold: the latter by Proposition~2.9 in~\cite{CaMob}, using again cocompactness of $H$.

We can now apply Theorem~7.4 in~\cite{CaMoa} and deduce that each $L_i$ acts cocompactly on $X_i$; indeed the proof
of \emph{loc. cit.} even provides a quasi-isometry bewteen $X_i$ and the model space (symmetric space or
Bruhat--Tits building) of $L_i$. Thus $L$ acts cocompactly on $X$, which implies that $L$ is cocompact in $H$.

Theorem~7.4 in~\cite{CaMoa} also provides a Tits-isometric identification of $\bd X_i$ with the boundary of
the model space of $L_i$. Assuming now that each $L_i$ has rank at least two, we can apply Tits' rigidity
theorem (Theorem~5.18.4 in~\cite{Tits74}) and deduce that $\Isom(X_i)$ is the group of isometries of the model
space, which coincides with the group of automorphisms of the associated semi-simple group.
\end{proof}

\section{Presence of an algebraic factor}
\subsection{Algebraic factors in general}\label{sec:algebraic}
Following Margulis~\cite[IX.1.8]{Margulis}, we shall say that a simple algebraic group $\mathbf{G}$ defined over
a field $k$ is \textbf{admissible} if \emph{none} of the following holds:
\begin{itemize}
\item[---] $\chr(k) = 2$ and $\mathbf{G}$ is of type $A_1$, $B_n$, $C_n$ or $F_4$;

\item[---] $\chr(k) = 3$ and $\mathbf{G}$ is of type $G_2$.
\end{itemize}
A semi-simple group will be called admissible if all its simple factors are.

\begin{thm}\label{thm:algebraic:factor}
Let $k$ be a local field and $\GG$ an adjoint admissible connected semi-simple $k$-group without
$k$-anisotropic factors. Let $X$ be a proper \cat space without Euclidean factor and $H<\Isom(X)$ a closed totally disconnected subgroup
acting minimally and without fixed point at infinity. Let $\Gamma<\GG(k)\times H$ be any lattice; in case $\rank_k\GG=1$ and
$\chr(k)>0$, we assume in addition $\Gamma$ cocompact.

If $\Gamma$ projects faithfully to the simple factors of $\GG(k)$, then $H$ is a semi-simple algebraic group upon passing to a finite
covolume subgroup containing the image of $\Gamma$. Furthermore, $\Gamma$ is finitely generated.
\end{thm}

\begin{remark}
There is a similar statement without the \cat space $X$ in Theorem~5.20 of~\cite{CaMob}, but at the cost of
assuming $H$ compactly generated and assuming that $\Gamma$ projects densely to $H$. In general, we do not know
how to prove \emph{a priori} that the closure of the projection of $\Gamma$ to $H$ is compactly generated, even
if we assume $H$ compactly generated. In the above theorem, compact generation of $H$ follows \emph{a
posteriori} from the statement. In fact, the bulk of the proof given below is concerned with addressing this very
issue.
\end{remark}

We begin with a geometric finiteness result that will allow us to \emph{rule out phenomena of ad\'elic type} in the proof
of Theorem~\ref{thm:algebraic:factor}; for its own sake, we provide more generality.

\begin{prop}\label{prop:NoAdele}
Let $X$ be a proper \cat space and $H<\Isom(X)$ a closed subgroup acting minimally and without fixed point at
infinity. Let $\{H_n\}$ be a non-decreasing family of closed subgroups of $H$ such that the closure of the union
of all $H_n$ is co-amenable in $H$.

Then there is $N\in\NN$ such that no $H_n$ can be a compact extension of a direct product of more than $N$
non-compact factors.
\end{prop}

\begin{proof}
Let $X = Y\times \RR^d$ be the maximal Euclidean decomposition, so that
$\Isom(X)=\Isom(Y)\times\big(\OO(d)\ltimes\RR^d\big)$, see Theorem~II.6.15 in~\cite{Bridson-Haefliger}. Arguing
by contradiction, we can assume that each $H_n$ has a compact normal subgroup $K_n$ such that $H_n/K_n$ can be
decomposed as a direct product of $n$ non-compact factors. We claim that we
can pass to a further subsequence and regroup factors so that all $n$ factors have non-compact image in $\Isom(Y^{K_n})$.
Indeed, each $H_n/K_n$ acts on a Euclidean subspace of $\RR^d$, namely its $K_n$-fixed points. This implies that at most $d$
of the non-compact factors of $H_n/K_n$ have a non-compact image in $(\RR^d)^{K_n}$; thus at least $n-d$ factors have non-compact
image in $\Isom(Y^{K_n})$, which implies the claim.

Since the closure $H_\infty<H$ of the union of all $H_n$ is co-amenable, it has no fixed point in $\bd Y$ by
Proposition~2.1 in~\cite{CaMob}. Therefore, by compactness of $\bd Y$, we can further assume that none of the
$H_n$ has a fixed point in $\bd Y$. It follows that each $H_n$ admits some minimal non-empty closed convex
invariant subspace $Y_n\se Y$ and that moreover the union $Z_n\se Y$ of all such subspaces splits isometrically
and equivariantly as $Z_n\cong Y_n\times T_n$, where the ``space of components'' $T_n$ is a \emph{bounded} \cat
space endowed with the trivial $H_n$-action; for all this, see Remark~39 in~\cite{Monod_superrigid}.

We claim that the sequence $\{T_n\}$ is of non-increasing diameter. Indeed, let $t, t'\in T_{n+1}$; then both
$Y_{n+1}\times\{t\}$ and $Y_{n+1}\times\{t'\}$ contain some, \emph{a priori} several, minimal $H_n$-subspaces.
We denote by $s, s'\in T_n$ the elements corresponding to some arbitrary such choices $Y_n\times\{s\}\se
Y_{n+1}\times\{t\}$ and $Y_n\times\{s'\}\se Y_{n+1}\times\{t'\}$. Now we have $d(t, t')\leq d(s, s')$ and the
claim follows.

In view of the claim, we may choose a sequence of points $y_n\in Y_n$ that remains bounded. Notice that $K_n$
acts trivially on $Y_n$. Our assumption on $H_n/K_n$ together with the splitting theorem
from~\cite{Monod_superrigid} shows that $Y_n$ admits a splitting as a product of $n$ non-compact factors. In
particular, we can choose $n$ geodesic rays issuing from $y_n$ and spanning a Euclidean $n$-dimensional
quadrant. Having Euclidean quadrants of unbounded dimension but based at the points $y_n$ which remain in a
bounded set contradicts the local compactness of $Y$.
\end{proof}

\begin{proof}[Proof of Theorem~\ref{thm:algebraic:factor}]
In view of the nature of the statement, we may and shall replace $H$ by the closure of the projection of
$\Gamma$, which has finite covolume in $H$.
By Theorem~2.4 in~\cite{CaMob}, $H$ still acts minimally and without fixed point at infinity.
In particular, we can assume it non-compact since otherwise it is
trivial, in which case the statement is empty except for the finite generation of $\Gamma$; the latter would still follow
as explained below for $\Gamma_U$, which coincides with $\Gamma$ when $H$ is trivial.

As we shall see, given the results we proved in~\cite{CaMob}, the main step here is to prove the following.

\smallskip
\emph{Main claim: the lattice $\Gamma$ is finitely generated.}

\smallskip

To this end, let $U<H$ be a compact open subgroup, which exists by~\cite[III.4 No~6]{BourbakiTGI}. Since
$\Gamma$ projects injectively to $\GG(k)$, we can consider $\Gamma_U=\Gamma\cap(\GG(k)\times U)$ as a lattice in $\GG(k)$.
Moreover, $\Gamma_U$ is irreducible in $\GG(k)$ since it projects injectively to the simple factors (recalling
that the various notions of irreducibility coincide in the case of lattices in semi-simple groups). Our assumptions
imply that $\Gamma_U$ is finitely generated; indeed, we recall the argument given in~\cite[5.20]{CaMob}: either
we have simultaneously $\chr(k)>0$ and $\GG$ is simple of $k$-rank one, in
which case we assumed $\Gamma$ cocompact, so that $\Gamma_U$ is cocompact in the compactly generated group
$\GG(k)$ and hence finitely generated~\cite[I.0.40]{Margulis}; otherwise, $\Gamma_U$ is known to be finitely
generated by applying, as the case may be, either Kazhdan's property, or the theory of fundamental domains, or
the cocompactness of $p$-adic lattices~--- we refer to Margulis, Sections~(3.1) and~(3.2) of Chapter~IX
in~\cite{Margulis}. 

We choose a non-decreasing sequence $\{\Gamma_n\}$ of finitely generated subgroups with
$\Gamma_U<\Gamma_n<\Gamma$ and which exhaust all of $\Gamma$. We denote by $G_n<\GG(k)$ the closure of the
projection of $\Gamma_n$ to $\GG(k)$ and by $H_n<H$ the closure of the projection of $\Gamma_n$ to $H$. Notice
that the closure of the union of all $H_n$ coincides with the closure of the projection of $\Gamma$
to $H$ and thus is all of $H$ in view of our preliminary reduction.

\smallskip

We claim that $\Gamma_n$ is a topologically irreducible lattice in $G_n\times H_n$ upon discarding the first few $n$.
Lemma~\ref{lem:intermediate:lattice} shows that $\Gamma_n$ is indeed a lattice in $G_n\times H_n$ and hence
the point to verify is that $G_n, H_n$ are both non-discrete.

If all $G_n$ are discrete, they are lattices in $\GG(k)$ and thus $\Gamma_U$ has finite index in $\Gamma_n$
since the projection of $\Gamma$ to $\GG(k)$ is faithful; in particular, $H_n$ is compact and hence fixes a
point in $X$. Considering the nested sequence of $H_n$-fixed points in the compactification $\overline{X}$, we
deduce by compactness that $H$ fixes a point in $\overline{X}$. This is impossible since $H$ fixes no point at infinity
and is non-compact.

If $H_n$ is discrete, then $\Gamma_n\cap (G_n\times 1)$ is a lattice (see Theorem~1.13 in~\cite{Raghunathan}).
Viewing it in $G_n$, it is a normal (hence cocompact) lattice since it is normalised by the projection of
$\Gamma_n$. However, $G_n$ does not admit a normal lattice when it is non-discrete. Indeed, being Zariski-dense
in $\GG$ (by Borel density applied to $\Gamma_U$), it contains the group $\GG_\alpha(k)^+$ for some simple
factor $\GG_\alpha$ by~\cite{Prasad77} and the latter is simple~\cite{Ti64} (and non-discrete). The
claim that $\Gamma_n$ is irreducible in $G_n \times H_n$ is proved.

\smallskip

We can now apply Theorem~5.1 from~\cite{CaMob} and deduce that $H_n$ is a compact extension of a semi-simple
algebraic group without compact factors. In fact, this reference allows a priori for a possibly infinite
discrete direct factor in $H_n$ which is also virtually a direct factor of $\Gamma_n$, but in the case at hand
this contradicts the fact that it is Zariski dense in a simple algebraic group, namely any simple factor $\GG_\alpha$
(since it contains $\Gamma_U$ which is Zariski-dense in $\GG$ by Borel density).

\smallskip

We claim that the obtained semi-simple quotient of $H_n$ is a direct factor of the quotient associated to $H_{n+1}$.

Indeed, Margulis' commensurator arthmeticity~\cite[1.(1)]{Margulis} shows that $\Gamma_U$ is $S$-arithmetic and
hence the projection of $\Gamma$ is contained in $\GG(K)$ for some global field $K$ over which $\GG$ is defined
by Theorem~3.b in~\cite{Borel66} (see also~\cite[Lemma~7.3]{Wortman}). An examination of the proof of
Theorem~5.1 in~\cite{CaMob} shows that the semi-simple quotient of $H_n$ is the product of the non-compact
semi-simple factors of all $\GG(K_v)$, where $v$ ranges over the set of places of $K$ for which $\Gamma_n$ is unbounded.
This proves the claim.

\smallskip

Proposition~\ref{prop:NoAdele} applies and we deduce from the previous claim that the sequence of the
semi-simple quotients associated to $\{H_n\}$ stabilises.
In view of the above discussion, it follows that $\Gamma<\GG(K)$ is in fact itself $S$-arithmetic.
In view of the assumptions on $\GG$ and of the results in Section~3.2 of Chapter~IX in~\cite{Margulis},
this $S$-arithmetic group is finitely generated if it is irreducible. Since $\Gamma$
projects injectively into the simple factors of $\GG(k)$, irreducibility
follows. (Alternatively, argue as in the proof of Proposition~\ref{prop:irred}.)
This concludes the proof of the main claim.

\smallskip

We now have $\Gamma=\Gamma_n$ for $n$ large enough; in particular, $H_n=H$ and the proof is complete.
\end{proof}

\subsection{Reduction to the totally disconnected case}

Retain the notation of \S\,\ref{sec:set-up}. The following result will later allow us to focus  on the case where
each $G_i$ is totally disconnected.

\begin{prop}\label{prop:Reduction}
Assume that the projection of $\Gamma$ to each $G_i$ is faithful.

If $G$ is not totally disconnected, then each $G_i$ contains a closed subgroup $H_i$ of finite covolume which is
a simple algebraic group over a local field and $\Gamma$ is $S$-arithmetic. If in addition $G$ acts cocompactly on
$X$, then $G_i/H_i$ is compact.
\end{prop}

\begin{proof}
Define $H_i$ as the closure of the projection of $\Gamma$ to $G_i$; we shall focus on the statements about $G_i$ and $H_i$,
since the arithmeticity of $\Gamma$ will then follow by Margulis' results
(see Theorem~(A) in Chapter~IX of~\cite{Margulis}).

If the identity component $G^\circ$ is non-trivial, then the same holds for some $G_i$. Upon renumbering the
$G_i$'s, we may and shall assume that $G_i$ is totally disconnected if and only if $i>k$ for some $k \in \{1,
\dots, n\}$. By~\cite[Corollary~1.11]{CaMoa}, it follows that $G_i$  is a non-compact simple Lie group with trivial
centre for each $i \in \{1, \dots, k\}$.

By Proposition~\ref{prop:irred}(i), the group $H_i$ is non-discrete for each $i$. In particular we have $H_i = G_i$ for each $i
\leq k$ by Borel density~\cite{Borel60}, since a Zariski-dense subgroup of a simple Lie group is either
discrete or dense. Furthermore, since $H_i$ has finite covolume in $G_i$, it follows from
\cite[Theorem~2.4]{CaMob} that $H_i$ acts minimally without fixed point at infinity on $X_i$. In particular it
has no non-trivial compact normal subgroup. Now Theorem~\ref{thm:algebraic:factor} implies that
$$H := G_1 \times \cdots \times G_k \times H_{k+1} \times \cdots \times H_n$$
is a semi-simple algebraic group. In view of Proposition~\ref{prop:semi-simple:cocompact}, if $G_i$ acts
cocompactly on $X_i$, then so does $H_i$.
\end{proof}

It will be convenient to have the following \emph{ad hoc} simpler variant of Proposition~\ref{prop:Reduction};
it is essentially just a shortcut available in positive characteristic.

\begin{prop}\label{prop:Reduction_ad_hoc}
Let $k$ be a local field of positive characteristic and $\GG$ an adjoint connected absolutely almost simple
$k$-group of positive $k$-rank. Let $X$ be a proper \cat space without Euclidean factor and let $H<\Isom(X)$ be a closed
subgroup acting cocompactly, minimally and without fixed point at infinity.

If there is any lattice $\Gamma<\GG(k)\times H$ that projects faithfully to $\GG(k)$, then $H$ is totally disconnected.
\end{prop}

\begin{proof}
Theorem~1.6 in~\cite{CaMoa} implies that $H$ is of the form
$H=S\times D$, where $S$ is a connected semi-simple Lie group and $D$ is totally disconnected. Let $U<D$ be a compact open subgroup
and observe that the lattice $\Gamma_U<\GG(k)\times S\times U$ (as considered in~\ref{sec:open}) still projects injectively to $\GG(k)$. 
Suppose for a contradiction that $S$ is non-compact. Then we have obtained
a lattice in $\GG(k)\times S$ which is irreducible and $S$-arithmetic in view of Margulis' arithmeticity~\cite{Margulis}.
This is absurd since the characteristics of the fields of definition do not agree.
\end{proof}

\subsection{Arithmeticity of residually finite lattices}\label{sec:residually_finite}
We remain in the setting of \S\,\ref{sec:set-up}.

\begin{thm}\label{thm:technical:ResFi}
Suppose that the lattice $\Gamma< G = G_1 \times \cdots \times G_n$ is algebraically irreducible. Assume that $G_1$
possesses an open subgroup $P$ which is a compact extension of a non-compact admissible simple algebraic group $\mathbf{H}$
over a local field $k$. In case $k$ has positive characteristic and $\HH$ has $k$-rank one, we assume in addition $\Gamma$ cocompact.

If $\Gamma$ is residually finite, then each $G_i$ contains a closed subgroup of finite covolume
which is a compact extension of a simple algebraic group over a local field, and $\Gamma$ is $S$-arithmetic.
\end{thm}

\begin{proof}
By Proposition~\ref{prop:ResFinite:irred}, the projection of $\Gamma$ to each $G_i$ is faithful. Therefore,
in view of Proposition~\ref{prop:Reduction}, we can assume that $G$ is totally disconnected.
Set
$$\Gamma_{P} = \Gamma \cap (P \times G_2 \times \cdots \times G_n).$$
By assumption $P$ has a compact normal subgroup $K$ such that $P/K=\HH(k)$. Now $\Gamma_{P}$ maps onto a lattice in the product
$$\HH(k) \times G_2 \times \dots \times G_n$$
and this map has finite kernel. Since $\Gamma_P$ is residually finite, we can assume that the kernel is trivial upon replacing
$\Gamma_P$ with a finite index subgroup; we henceforth consider $\Gamma_P$ as a lattice in the above product.

The projection of $\Gamma_P$ to $\HH(k)$ is faithful since we have already recorded that $\Gamma$ projects injectively to $G_1$.
Therefore, we can apply a first time Theorem~\ref{thm:algebraic:factor} to $\Gamma_P$ and deduce in particular for each $i\geq 2$
that $G_i$ is an admissible semi-simple algebraic group upon replacing it by a closed subgroup of finite covolume containing the image of $\Gamma_P$.
In fact, these groups are simple in view of the irreducibility of $X_i$ (\emph{e.g.} by the splitting theorem). We write $G_i=\GG_i(k_i)$ for
$i\geq 2$ and also note that $\Gamma_P$ is irreducible (\emph{e.g.} by Proposition~\ref{prop:irred}).

\smallskip

We now return to the lattice $\Gamma<G$ with the intention to apply a second time Theorem~\ref{thm:algebraic:factor}, but reversing the
r\^oles of $G_1$ and $G_2\times\cdots\times G_n$. We point out that the simple groups $\GG_i$ are all admissible since
both the absolute type and characteristic are constant over all factors in view of the fact that the $S$-arithmetic group $\Gamma_P$ is
irreducible.
However, we have no guarantee that the technical assumption made on $\HH$ holds for $\GG_i$. It can indeed fail and likewise the finite generation
used in the proof of Theorem~\ref{thm:algebraic:factor} for $\Gamma_U$ is also known to fail. We shall now circumvent this difficulty.

The group $\Gamma_P$ is finitely generated by the above application of Theorem~\ref{thm:algebraic:factor}. We consider a non-decreasing sequence of
finitely generated groups $\Gamma_j$ starting with $\Gamma_0=\Gamma_P$ and exhausting $\Gamma$. We denote by $L_j$ and $R_j$ the closure of
the projection of $\Gamma_j$ to
$$G_1\times \cdots \times G_{n-1} \kern5mm\text{and}\kern5mm  G_n$$
respectively. Lemma~\ref{lem:intermediate:lattice} shows that $\Gamma_j$ is a lattice in $L_j\times R_j$.
It is topologically irreducible in the former product
since already the projections of $\Gamma_0$ are non-discrete (\emph{e.g.} by Proposition~\ref{prop:NonDiscreteProj}).
Theorem~5.1 in~\cite{CaMob} implies that $L_j$ is a compact extension of a semi-simple group $S_j$. We write $Q=P\times G_2\times\cdots\times G_{n-1}$,
wherein $Q=P$ is understood if $n=2$. The group $Q\cap L_j$ is open in $L_j$ and non-compact since it contains $L_0$ which is of finite covolume
in the non-compact group $Q$. Therefore the image of $Q\cap L_j$
in $S_j$ contains $S_j^+$ by~\cite[thm.~(T)]{Prasad82} (or by an application of Howe--Moore). Since $S_j^+$ is cocompact in $S_j$
(see~\cite[6.14]{Borel-Tits73}) we conclude that $Q\cap L_j$ is cocompact in $L_j$ and hence has finite index. This shows that $L_0$ has
finite index in $L_j$. Therefore, $S_0$ has finite index in $S_j$ for all $j$; it follows, since $S_0^+$ is
simple~\cite{Ti64}, that $S_j$ normalises $S_0^+$. We denote by $L_0^+$ the preimage of $S_0^+$ in $L_0$ and set
$\Gamma_0^+=\Gamma\cap (L_0^+\times G_n)$. Then $\Gamma_0^+$ has finite index in $\Gamma_0$ since the latter is finitely generated and
since $S_0/S_0^+$ is a virtually Abelian torsion group~\cite[6.14]{Borel-Tits73}. Moreover, $\Gamma_j$ normalises $\Gamma_0^+$ for all $j$
in view of corresponding statement for $S_0^+<S_j$ above.

At this point we have a natural map $\Gamma\to \Aut(S_0^+)$ whose image normalises $\Gamma_0^+$. In combination with the injective map
$\Gamma\to G_n$, we have a realisation of $\Gamma$ in the normaliser of $\Gamma_0^+$ in the algebraic group $\Aut(S_0^+)\times \GG_n(k_n)$.
Since $\Gamma_0^+$ is a lattice in the latter, it follows from Borel's density theorem that $\Gamma_0^+$ has finite index in its normaliser,
see~\cite[II.6.3]{Margulis}. In conclusion, $\Gamma_0^+$ has finite index in $\Gamma$ and thus both sequences $\Gamma_j$ and $L_j$ are eventually
constant, completing the proof.
\end{proof}

\section{Kac--Moody groups}
\subsection{A lemma on Coxeter--Dynkin diagrams}
A theorem of G.~Moussong characterises the Gromov hyperbolic Coxeter groups in terms of their Coxeter diagram.
In fact, Moussong's result says that a finitely generated Coxeter group is Gromov hyperbolic if and only if it
does not contain any subgroup isomorphic to $\ZZ \times \ZZ$. The latter property can easily be detected on the
Coxeter diagram of $G$, since any subgroup isomorphic to $\ZZ \times \ZZ$ is conjugate into a special parabolic
subgroup of $W$ which is either of affine type or which decomposes as the direct product of two infinite special
subgroups. It turns out that for some specific families of Coxeter groups, the presence of a $\ZZ \times
\ZZ$-subgroup always guarantees the presence of an affine parabolic.

\begin{lem}\label{lem:Dynkin}
Let $(W, S)$ be a crystallographic Coxeter system of simply laced or $3$-spherical type, with $S$ finite. Then
$W$ is Gromov hyperbolic if and only if $W$ contains no parabolic subgroup of affine type.
\end{lem}

Recall that a Coxeter group $W$ is called \textbf{crystallographic} if its natural geometric representation in
$\RR^n$  preserves a lattice (see~\cite{Bourbaki_Lie456}). This property is known to be equivalent to each of
the following conditions:
\begin{itemize}
\item The Coxeter numbers which appear in a Coxeter presentation of $W$ belong to $\{2, 3, 4, 6, \infty\}$.

\item $W$ is the Weyl group of some Kac--Moody Lie algebra.
\end{itemize}

In particular, if $W$ is the Weyl group of a Kac--Moody group over any field, then $W$ is crystallographic.

\begin{proof}[Proof of Lemma~\ref{lem:Dynkin}]
We may assume  that $W$ is irreducible. If $W$ possesses a parabolic subgroup of affine type, then it contains a
$\ZZ \times \ZZ$-subgroup and cannot be hyperbolic. Assume now that $W$ is not hyperbolic. In view of Moussong's
theorem, all we need to show is that if $W$ contains two infinite special subgroups $W_{I}, W_{J}$ which
mutually commute, then it also contains a parabolic subgroup of affine type. Without loss of generality we may
assume that $W_I$ and $W_J$ are \emph{minimal} infinite special subgroups, namely that every proper special
subgroup of $W_I$ or $W_J$ is finite. The list of minimal infinite Coxeter groups is known and may be found in
Exercises~13--17 for \S\,4 in Chapter~V from~\cite{Bourbaki_Lie456}. It turns out that every such a Coxeter group
is either affine or is defined by a diagram belonging to a short list, the members of which have size~$\leq 5$.
A short glimpse at this list shows that none of them is simply laced. Furthermore, only three of them are
$3$-spherical crystallographic, namely those depicted in Figure~\ref{fig:Dynkin}.

\begin{figure}[h]
\begin{center}
\begin{tabular}{c @{\hspace{12mm}}c @{\hspace{12mm}} c }
\includegraphics[width=25mm,height=25mm]{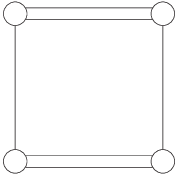}
&

\includegraphics[width=25mm,height=25mm]{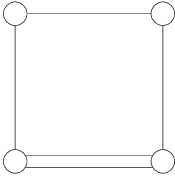}

&

\includegraphics[width=37mm,height=33mm]{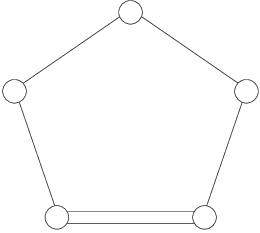}\\

(a) & (b) & (c)
\end{tabular}
\end{center}
\caption{Minimal non-spherical $3$-spherical Dynkin diagrams}%
\label{fig:Dynkin}
\end{figure}

We now consider a path of minimal possible length joining $I$ to $J$ in the Coxeter diagram of $W$ and consider
the Coxeter diagram induced on the union of the vertex set $P $ of this path together with $I \cup J$.

If any vertex in $P -(I \cup J)$ is linked by an edge to a vertex of $I$ or $J$ belonging to an edge
labelled~$4$, then the Coxeter diagram contains a subdiagram of type $\tilde B_2$ or $\tilde C_2$ and we are
done. Similarly, if some vertex of $P -(I \cup J)$ is linked to a vertex of $I$ or $J$ by an edge labelled~$4$,
then we are done as well. Thus we may assume that all labels of edges linking a vertex in $P -(I \cup J)$ to a
vertex in $I \cup J$ are~$3$, and that all such edges are not adjacent to an edge labelled~$4$ in $I$ or $J$.

Now it follows that if more than one vertex of $I$ or $J$ is linked by an edge to vertex in $P -(I \cup J)$,
then the diagram contains a subdiagram of type~$\tilde A_n$ and we are done. It remains to consider the case
where each vertex of $P -(I \cup J)$ is linked to at most one vertex in $I$ and $J$. In that case, it is readily
seen that the diagram contains a subdiagram of type~$\tilde C_n$. This finishes the proof.
\end{proof}

Another useful and well-known fact is the following.

\begin{lem}\label{lem:Dynkin:2}
Let $(W, S)$ be an irreducible non-spherical non-affine Coxeter system such that for each proper subset $J
\subset S$, the special subgroup $W_J$ is either spherical or affine. Then $|S| \leq 10$.
\end{lem}

\begin{proof}
See Exercises~13--17 for \S\,4 in Chapter~V from~\cite{Bourbaki_Lie456}.
\end{proof}

\subsection{Complete Kac--Moody groups and their buildings}\label{sec:KM}
Let $G$ be a complete  adjoint Kac--Moody group over a finite field $\FF_q$. Such a group may be obtained as
follows. Start with a Kac--Moody--Tits functor $\mathcal{G}$ associated to a Kac--Moody root datum of adjoint
type, as defined in~\cite{Tits87} (see also~\cite{RemAst} for non-split versions). Thus $\mathcal{G}$ is a group
functor on the category of commutative rings. Its value on any field $k$ is a group $\mathcal{G}(k)$ which
admits two natural uniform structures. Completing $\mathcal{G}(k)$ with respect to any of these yields a totally
disconnected topological group $\mathbf{G}(k)$ which contains $\mathcal{G}(k)$ as a dense subgroup, see
\cite[\S\,1.2]{CaRe}. When $k= \FF_q$ is a finite field of order $q$, then $\mathbf{G}(k)$ is locally compact.

We remark that the functor $\mathbf{G}$ may be obtained by a different construction, due to Olivier
Mathieu~\cite{Mathieu}, which yields not only a group functor but an ind-group scheme.

\medskip%
We assume henceforth that $k = \FF_q$ and set $G = \mathbf{G}(\FF_q)$. The group $G$ possesses a $BN$-pair with
$B$ compact open in $G$. This $BN$-pair yields a locally finite building $X$ of thickness $q+1$ on which $G$
acts faithfully, continuously and properly by automorphisms. Furthermore $X$ possesses a natural realisation as
a \cat space whose isometry group contains $\Aut(X)$ as a closed subgroup~\cite{Davis}. By a slight abuse of
notation, we shall not distinguish between $X$ and its \cat realisation. Thus we view $G$ as a closed subgroup
of $\Isom(X)$. The $G$-action on $X$ is transitive on the chambers and the chambers are compact. In particular
$G$ stabilises a minimal closed convex invariant nonempty subspace which we may view as a \cat realisation of
$X$ on which $G$ acts minimally. There is thus no loss of generality in assuming that $G$ acts minimally on $X$.

The fact $G$ has no fixed point at infinity may be established in several different ways. The conceptually
easiest one is the following. It is shown in~\cite[Lemma~9]{CaRe} that the derived group $[G, G]$ is dense in
$G$. It follows that if $G$ fixed a point $\xi$ in the boundary at infinity $\bd X$, then $G$ would stabilise
each horoball centered at that point, contradicting the minimality of the action. Another way to obtain this
statement is by using the fact that a Coxeter group has no fixed point at infinity in its natural action on the
associated \cat cell complex as one sees by considering the numerous reflections (or else by applying
\cite[Theorem~3.14]{CaMob}). Since for any apartment $A \subset X$ the $\Stab_G(A)$-action on $A$ is isomorphic
to the natural action of the Weyl group on its cell complex, and since apartments are convex, it follows again
that $G$ has no fixed point at infinity.

Finally, since $X$ has a cocompact isometry group it has finite-dimensional Tits boundary by
\cite[Theorem~C]{Kleiner}. This discussion shows in particular that the group $G_1$ appearing in the statement of Theorem~\ref{thm:RankRigid} satisfies
the set-up described in \S\,\ref{sec:set-up} by considering its natural action on the associated building.

\medskip
Furthermore, it turns out that $G$ is topologically simple~\cite[Proposition~11]{CaRe}. In
addition, if the ground field $\FF_{q}$ has order $q \geq 1764^{d}/25$, where $d$ denotes the
dimension of the building $X$, and if $W$ is $2$-spherical, then $G$ has Kazhdan's property~(T)
\cite[Corollary~G]{DJ02}. Notice that the dimension of $X$ is bounded above by the maximal rank of a finite
Coxeter subgroup of $W$, see~\cite{Davis}.

\begin{lem}\label{lem:parabolic}
Let $G$ be an irreducible complete Kac--Moody group of adjoint type over a finite field $\FF_{q}$. Assume that
the Weyl group $W$ of $G$ is infinite and simply laced or $3$-spherical but not Gromov hyperbolic. Then $G$ contains an open
subgroup $P$ which is a compact extension of a simple algebraic group over a local field of characteristic~$p =
\charact \FF_q$ and rank~$\geq 2$. Furthermore, if $W$ is simply laced or if $\charact \FF_q \neq 2$, then the
latter simple group is admissible.
\end{lem}

\begin{proof}
By Lemma~\ref{lem:Dynkin}, the group Weyl group $W$ possesses a special parabolic subgroup of (irreducible)
affine type $W_J$. Let $P_J < G$ be a parabolic subgroup of type $W_J$. Thus $P_J = B \cdot W_J \cdot B$, where
$B$ denotes the Borel subgroup of $G$, namely the $B$-subgroup of the $BN$-pair. In particular $P_J$ contains
the compact open subgroup $B$ and is thus open. The subgroup
$$K_J = \bigcap_{g \in P_J} g B g\inv$$
is a compact normal subgroup of $P_J$. The quotient $P_J/K_J$ is a complete Kac--Moody group of type $W_J$ over
$\FF_q$ (see~\cite[Proposition~11]{CaRe} and~\cite[\S\,5]{CER}). It follows from~\cite[Appendix]{Tits84} (or else
from the uniqueness theorem in~\cite{Tits87}) that $P_J/K_J$ is a simple algebraic group whose Weyl group is the
spherical Weyl group of $W_J$. This yields the desired conclusions.
\end{proof}

\begin{remark}\label{rem:BorelHarder}
As pointed out by G.~Margulis~\cite[IX.1.6(viii)]{Margulis}, it follows from the arithmeticity theorem, combined
with~\cite[Korollar~1 p.~133]{Harder}, that if $W_J$ is not of type $\tilde A_n$, then $P_J$ does \emph{not} admit
any uniform lattice. (For type $\tilde A_n$, such lattices indeed exist, see~\cite{BorelHarder} and
\cite{Cartwright-Steger}.) If follows in particular that if $W_J$ is not  of type $\tilde A_n$, then \emph{no}
product of the form $G \times H$, where $H$ is a totally disconnected locally compact group, possesses
\emph{any} uniform lattice.

Notice furthermore that the condition that every special subgroup of $W$ of affine type  be of type $\tilde A_n$
is a very strong one. For example, if $W$ is $3$-spherical, then only finitely many Coxeter diagrams are
possible for $W$. This shows that in general one should not expect $G$ (or $G \times H$) to possess any uniform
lattice.
\end{remark}

\section{Completion of the proofs}
\subsection{Reduction of the hypotheses}
For the purposes of this last section, let us define a complete Kac--Moody group $G = \mathbf{G}(\FF_q)$ over a finite
field $\FF_q$ with Weyl group $W$ to be \textbf{admissible} if any of the following two conditions holds:
\begin{itemize}
\item[---] $W$ is simply laced.

\item[---]  $\charact(\FF_q) \neq 2$ and $W$ is $3$-spherical but not Gromov hyperbolic.
\end{itemize}

Notice that the Weyl group $W$ of $\mathbf{G}$ is Gromov hyperbolic if and only if $\mathbf{G}(F)$ is Gromov
hyperbolic for each finite field $F$. Indeed $\mathbf{G}(F)$ acts properly and cocompactly on a building of type
$W$, and it is known that a building is Gromov hyperbolic if and only if its Weyl group is so (see \emph{e.g.}~\cite{Davis}).

Although we have already used the term admissible in a different context in \S\,\ref{sec:algebraic}, the above
definition will cause no confusion. Indeed, given a Kac--Moody group $\mathbf{G}$ of affine type over $\FF_q$
(equivalently $\mathbf{G}(\FF_q)$ is isomorphic to a semi-simple algebraic group $\mathbf{H}$ over a field $k$
of formal power series with coefficients in $\FF_q$), if $\mathbf{G}(\FF_q)$ is admissible in the above sense
then $\mathbf{H}(k)$ is admissible in the sense of \S\,\ref{sec:algebraic}.

\bigskip

We now proceed to relate the broad assumptions of the Introduction to the setting considered in \S\,\ref{sec:set-up}.

Let $n\geq 2$; for each $i \in \{1, \dots, n\}$, let $X_i$ be a proper \cat space and $G_i < \Isom(X_i)$ be a
closed subgroup acting cocompactly. We recall that cocompactness implies that $G_i$ is compactly generated
(see \emph{e.g.} Lemma~22 in~\cite{Mineyev-Monod-Shalom}).
Assume that $G_1$ is an admissible irreducible Kac--Moody group as discussed above.
Set $G=G_1\times\cdots\times G_n$ and $X=X_1\times\cdots\times X_n$. Finally, 
let $\Gamma<G$ be a lattice whose projection to each $G_i$ is faithful.
We assume $G_1$ infinite; this hypothesis was not made in Theorem~\ref{thm:RankRigid}
but the latter is otherwise trivial since $\Gamma$ would be finite.

\begin{prop}\label{prop:triv:euclidean}
The space $X$ has trivial Euclidean factor and $G$ has no fixed point at infinity.

Moreover, for each $i \in \{1, \dots, n\}$, there is a closed normal subgroup of finite index $G_i^* \unlhd G_i$, a
proper \cat space $Y_i = Y_{i, 1} \times \cdots \times Y_{i, k_i}$, where each $Y_{i, j}$ is
irreducible~$\neq\RR$ with finite-dimensional boundary and a continuous proper map $G_i^*
\to \Isom(Y_{i, 1}) \times \cdots \times \Isom(Y_{i, k_i})$ which yields a cocompact
minimal $G_i^* $-action on $Y_i$ without fixed point at infinity. Finally, for all $i,j$ the image of $G_i^*$ in $\Isom(Y_{i,j})$
is either totally disconnected or a connected simple Lie group.
\end{prop}

\begin{proof}
Since the $G_i$-action on $X_i$ is cocompact, there is a non-empty closed convex $G_i$-invariant subset $Y_i
\subseteq X_i$ on which the induced $G_i$-action is minimal. This action is proper and remains cocompact, 
which implies that the boundary $\bd Y_i$ is finite-dimensional (Theorem~C in~\cite{Kleiner}).
Corollary~5.3(ii) in~\cite{CaMoa} now states that $Y_i$ possesses a decomposition $Y_i = \RR^{d_i} \times Y_{i,
1} \times \cdots \times Y_{i, k_i}$, where $Y_{i, j}$ is an irreducible proper \cat space, such that
$$\Isom(Y_i) = \Isom(\RR^{d_i}) \times  \bigg( \big( \prod_{j=1}^{k_i} \Isom(Y_{i, j})\big) \rtimes F\bigg) ,$$
where $F$ is a finite permutation group of possibly isometric factors. Thus $G_i$ possesses a closed normal
subgroup of finite index $G_i^*=\prod_{j=1}^{k_i} G_{i,j}$ whose induced action on $Y_i$ is componentwise.

We now proceed to prove that $Y:=Y_1\times\cdots\times Y_n$ has no Euclidean factor, \emph{i.e.} $d_i=0$ for all $i$.
Our assumption on $G_1$ implies $d_1=0$ (see~\cite{Caprace-Haglund}). Considering the canonical Euclidean decomposition
of $Y$ (see~\cite[II.6.15]{Bridson-Haefliger}), we write $Y\cong Y'\times \RR^d$, where $Y'$ has no Euclidean factor and
$d=d_1+\cdots+d_n$. We claim that all $G^*$-fixed points at infinity lie in $\bd \RR^d$, where $G^*=\prod G_i^*$.
To this end, we observe that $\Gamma$ provides us with a lattice in $G^*$ upon passing to a finite index subgroup; we still
denote it by $\Gamma$. If $\Gamma$
is finitely generated, the claim follows from Proposition~3.15 in~\cite{CaMob}; in general, it is a consequence of the unimodularity
of $G^*$, a fact we establish in~\cite{Caprace-Monod_amen}.

We can now apply Theorem~1.6 from~\cite{CaMoa} and deduce that each $G_{i,j}$ is either totally disconnected or a connected
simple Lie group (modulo the compact kernel of its action on $Y_{i,j}$).
Proposition~3.6 in~\cite{CaMob} states that when $\Gamma$ is finitely generated, it virtually splits off
an Abelian direct factor of $\QQ$-rank $d$. The finite generation, however, is only used to provide a complementary factor
to this Abelian subgroup; the existence of a normal Abelian subgroup $A\lhd \Gamma$ of $\QQ$-rank $d$ is established in general
in \emph{loc.\ cit}.
We finally deduce that $d=0$ from the fact that $\Gamma$ projects injectively to the Kac--Moody group $G_1$ using~\cite{Caprace-Haglund}.

At this point we have established that $Y$ has trivial euclidean factor and that $G^*$ has no fixed points at infinity.
In particular $G$ has no fixed points in $\bd X=\bd Y$ and $X$ has no Euclidean factor either since $Y$ has finite codiameter in $X$.
\end{proof}

\subsection{End of the proofs}

\begin{proof}[Proof of Theorem~\ref{thm:alternative:1}]
Retain the notation of the theorem. Then Lemma~\ref{lem:parabolic} ensures that $G_1$ possesses an
open subgroup $P$ which is a compact extension of an admissible simple algebraic group of rank~$\geq 2$ over a local field.
We can assume $\Gamma$ residually finite. The statement of Theorem~\ref{thm:alternative:1} is not affected
by the reductions of Proposition~\ref{prop:triv:euclidean}; therefore, Theorem~\ref{thm:technical:ResFi} yields the
desired conclusion.
\end{proof}

\begin{proof}[Proof of Theorem~\ref{thm:RankRigid}]
We adopt the notation of the theorem.
Countrary to Theorem~\ref{thm:alternative:1}, the irreducibility assumption is in the present case sensitive to
replacing $G_i$ with the subfactors $G_{i,j}$ of Proposition~\ref{prop:triv:euclidean}. However, since $G_1$ is irreducible,
Proposition~\ref{prop:irred} implies that $\Gamma$ is at least algebraically irreducible.

As above, Lemma~\ref{lem:parabolic} provides an open subgroup $P<G_1$ which is a compact extension of an admissible simple algebraic group
$\HH(k)$ of rank~$\geq 2$ over a local field $k$; we emphasise that $k$ has positive characteristic.

The canonical image of $\Gamma_P=\Gamma\cap (P\times G_2\times \cdots\times G_n)$ in $\HH(k) \times G_2\times \cdots\times G_n$
is a lattice and it projects injectively to $\HH(k)$ in
view of the corresponding assumption on $\Gamma$. Thus Proposition~\ref{prop:Reduction_ad_hoc} implies that all $G_i$ are totally
disconnected. In particular $G_n$ possesses a compact open subgroup $U$ (see~\cite[III.4 No~6]{BourbakiTGI}). Set
$$\Gamma_{U} = \Gamma \cap (G_1 \times \cdots \times G_{n-1} \times U).$$
By assumption, the projection of $\Gamma_{U}$ to $U$ is faithful. In particular $\Gamma_{U}$ is residually
finite since $U$ is so, being a profinite group. Applying the faithfulness assumption to any other factor $G_i$,
we further deduce that  $\Gamma_{U}$ intersects the compact group $1 \times \cdots \times 1 \times U$ trivially;
therefore, we can view $\Gamma_U$ as a lattice in  the product
$$G_1 \times \prod_{j=1}^{k_2} G_{2,j}\times \cdots \times \prod_{j=1}^{k_{n-1}} G_{n-1,j}.$$
By Proposition~\ref{prop:CompactOpen:irred}, the group $\Gamma_{U}$ is algebraically irreducible. Thus we can apply
Theorem~\ref{thm:technical:ResFi} and deduce that $G_1$ and each factor $G_{i,j}$
contains a closed subgroup of finite covolume which is a simple algebraic group over a local field.

We now return to the initial lattice $\Gamma$ in $G$, which is algebraically irreducible by Proposition~\ref{prop:irred},
and conclude the proof as in the end of the proof of Theorem~\ref{thm:technical:ResFi}.
\end{proof}

\begin{proof}[Proof of Corollary~\ref{cor:simple}]
Since $\Gamma$ is irreducible and each $G_i$ is topologically simple (as recalled in Section~\ref{sec:KM}), it follows that the projection of $\Gamma$
to each $G_i$ is faithful. In view of Theorem~\ref{thm:RankRigid}, we may assume that $n=2$ and that $\Gamma$ is
not residually finite. All we need to show is that $\Gamma$ is virtually simple.

Since $G_1$ and $G_2$ are topologically simple and $\Gamma$ is irreducible, it follows from
\cite[Theorem~1.1]{Bader-Shalom} that if $\Gamma$ is uniform, then every non-trivial normal subgroup of $\Gamma$
has finite index. If $\Gamma$ is not uniform, then it has property~(T) in view of our assumptions and, hence,
the same conclusion on normal subgroups holds in view of~\cite[Theorem~1.3]{Bader-Shalom}.

Therefore Proposition~1 from~\cite{Wilson71} ensures that $\Gamma$ is virtually isomorphic to a direct product
of finitely many isomorphic simple groups. Since $\Gamma$ is irreducible as an abstract group by
Proposition~\ref{prop:irred}, we deduce that the latter direct product has a single simple factor. Thus $\Gamma$
is virtually simple.
\end{proof}

\begin{remark}
As pointed out by the anonymous referee, the above arguments show also the following. Let $\Gamma$ be a finitely generated group without non-trivial infinite index normal subgroup. Suppose that $\Gamma$ acts by isometries faithfully, minimally and without fixed point at infinity on an irreducible proper \cat space $X$. Then $\Gamma$ is either residually finite or virtually simple.

Indeed, in view of the above quoted result of Wilson, it suffices to show that $\Gamma$ does not have a finite index subgroup $\Gamma^*\cong \Gamma_1\times \Gamma_2$ splitting as a product of two infinite groups $\Gamma_i$. Since $X$ is irreducible, we can assume that it has no Euclidean factor for otherwise $X=\RR$ in which case the statement is obvious. Therefore our ``Borel density'' in the generality of Proposition~\ref{prop:LatticeConsequences} (presently invoked with $n=1$) implies that $\Gamma^*$ still acts minimally and without fixed point at infinity. By the splitting theorem of~\cite{Monod_superrigid}, this forces at least one of the $\Gamma_i$ to act trivially, a contradiction.

New examples of groups to which the above applies are provided in unpublished work of Shalom--Steger.
\end{remark}

\begin{proof}[Proof of Corollary~\ref{cor:MaximalLattices}]
In view of Corollary~\ref{cor:simple}, the assumption that the $G_i$'s are of non-affine type implies that $n=2$
in the above, and that any irreducible lattice $\Gamma$ of $G$ is virtually simple.  The finite residual
$\Gamma^{(\infty)}$ of $\Gamma$ is thus a normal subgroup of finite index, and any subgroup of $G$
commensurating $\Gamma$ normalizes $\Gamma^{(\infty)}$. Thus $\Comm_G(\Gamma) = \norma_G(\Gamma^{(\infty)})$.

Under the present hypotheses, the group $G$, and hence also $\Gamma$ has Kazhdan's property (T). Thus $\Gamma$
is finitely generated and, hence so is the lattice $\Gamma^{(\infty)}$. In view of~\cite[Corollary~2.7]{CaMob},
it follows that $\norma_G(\Gamma^{(\infty)})$ is itself a lattice in $G$, which is thus the desired maximal
lattice.
\end{proof}

\subsection{A lattice in a product of a simple algebraic group and a Kac--Moody group}\label{sec:example}
Let $G$ be an irreducible complete Kac--Moody group of simply laced type over a finite field $\FF_q$. It is
shown in~\cite{RemCRAS} (see also~\cite{CarGar} and~\cite{CaRe}) that the group $G \times G$ contains an
irreducible non-uniform lattice $\Gamma$, provided that $q$ is larger than the rank~$r$ of the Weyl group of
$G$. Assume now that $G$ is not of affine type. By Lemma~\ref{lem:parabolic} $G$ contains an open subgroup $P <
G$ which possesses a compact normal subgroup $K$ such that $P/K$ is a simple algebraic group over a local field. As in
the proof of Theorem~\ref{thm:technical:ResFi}, we may consider the group
$$\Gamma_P = \Gamma \cap (P \times G)$$
and view it as a lattice in $P/K \times G$. Furthermore $\Gamma_P$ is irreducible (see
Proposition~\ref{prop:CompactOpen:irred}) and one shows, as in the proof of
Theorem~\ref{thm:technical:ResFi}, that $\Gamma_P$ is finitely generated provided the ground field is large
enough.

Since $P/K$ acts cocompactly on the associated Bruhat--Tits building, we see in particular that $\Gamma_P$ is an
example of a \cat lattice (in the sense of~\cite{CaMoCRAS,CaMob}) in a product of an affine and a non-affine building.

Proposition~\ref{prop:ResFinite:irred} together with Theorem~\ref{thm:technical:ResFi} imply that
$\Gamma_P$ is \emph{not} residually finite. In particular its projection to the linear group $P/K$ is \emph{not}
faithful since a finitely generated linear group is residually finite by~\cite{Malcev40}. This example shows
that the assumption on the faithfulness of the projections in Theorem~\ref{thm:RankRigid} may not be removed.

%
%
\def\cprime{$'$}
\providecommand{\bysame}{\leavevmode\hbox to3em{\hrulefill}\thinspace}
\providecommand{\MR}{\relax\ifhmode\unskip\space\fi MR }
\providecommand{\MRhref}[2]{%
  \href{http://www.ams.org/mathscinet-getitem?mr=#1}{#2}
}
\providecommand{\href}[2]{#2}

\end{document}